\documentstyle[12pt]{article}
\begin{document}

\title{Isomorphisms of Hilbert C*-Modules and $*$-Isomorphisms of Related
Operator C*-Algebras}
\author{Michael Frank}
\maketitle

\begin{abstract}
\noindent
Let $\cal M$ be a Banach C*-module over a C*-algebra $A$ carrying
two $A$-valued inner products $\langle .,. \rangle_1$, $\langle
.,. \rangle_2$ which induce equivalent to the given one norms on $\cal M$.
Then the appropriate unital C*-algebras of adjointable bounded $A$-linear
operators on the Hilbert $A$-modules $\{ {\cal M}, \langle .,. \rangle_1 \}$
and $\{ {\cal M}, \langle .,. \rangle_2 \}$ are shown to be $*$-isomorphic if
and only if there exists a bounded $A$-linear isomorphism $S$ of these two
Hilbert $A$-modules satisfying the identity $\langle .,. \rangle_2 \equiv
\langle S(.),S(.) \rangle_1$. This result
extends other equivalent descriptions due to L.~G.~Brown, H.~Lin and
E.~C.~Lance. An example of two non-isomorphic Hilbert C*-modules with
$*$-isomorphic C*-algebras of ''compact''/adjointable bounded
module operators is indicated.
\end{abstract}

\noindent
Investigations in operator and C*-theory make often use of C*-modules as a
tool for proving, especially of Banach and Hilbert C*-modules.
Impressing examples of such applications are G.~G.~Kasparov's approach to
K- and KK-theory of C*-algebras \cite{JT,NEWO} or the investigations of
M.~Baillet, Y.~Denizeau and J.-F.~Havet \cite{BDH} and of Y.~Watatani
\cite{Wata} on (normal) conditional expectations of finite index on
W*-algebras and C*-algebras.
In addition, the theory of Hilbert C*-modules is interesting in its own.

\noindent
Our standard sources of reference to Hilbert C*-module theory are the papers
\cite{Pa1,Kas,Fr1,Frank:93}, chapters in \cite{JT,NEWO} and the book of
E.~C.~Lance \cite{Lance:95}.
We make the convention that all C*-modules of the present paper are
left modules by definition. A {\it pre-Hilbert $A$-module over a
C*-algebra} $A$ is an $A$-module $\cal M$ equipped with
an $A$-valued mapping $\langle .,. \rangle : {\cal M} \times
{\cal M} \rightarrow A$ which is $A$-linear in the first argument
and has the properties:
\[
\langle x,y \rangle = \langle y,x \rangle^* \; \, , \: \;
\langle x,x \rangle \geq 0 \quad {\rm with} \: {\rm equality} \: {\rm iff}
\quad x=0 \, .
\]
The mapping $\langle .,. \rangle$ is
called {\it the $A$-valued inner product on} $\cal M$. A pre-Hilbert
$A$-module $\{ \cal M, \langle .,. \rangle \}$ is {\it Hilbert} if and
only if it is complete with respect to the norm $\| . \| = \| \langle .,.
\rangle \|^{1/2}_A $. We always assume that the linear structures of
$A$ and $\cal M$ are compatible.

\noindent
One of the key problems of Hilbert C*-module theory is the question of
isomorphism of Hilbert C*-modules. First of all, they can be isomorphic as
Banach $A$-modules. But there is another natural definition: Two Hilbert
$A$-modules $\{ {\cal M}_1 , \langle .,. \rangle_1 \}$, $\{ {\cal M}_2 ,
\langle .,. \rangle_2 \}$ over a fixed C*-algebra $A$ are {\it isomorphic
as Hilbert C*-modules} if and only if there exists a bijective bounded
$A$-linear mapping $S: {\cal M}_1 \rightarrow {\cal M}_2$ such that the
identity $\langle .,. \rangle_1 \equiv$ $\equiv \langle S(.),S(.) \rangle_2$
is valid on ${\cal M}_1 \times {\cal M}_1$. In 1985 L.~G.~Brown presented two
examples of \hspace*{10cm}

\pagebreak[4]

\noindent
Hilbert C*-modules which are isomorphic as Banach C*-modules but which are
non-isomorphic as Hilbert C*-modules, cf.~\cite{Brown,Lin:90,Frank:93}.
This result was very surprising since Hilbert space theory, the classical
investigations on Hilbert C*-modules like \cite{Pa1,Kas},
\linebreak[4]
G.~G.~Kasparov's
approach to KK-theory of C*-algebras relying on countably generated Hilbert
C*-modules and other well-known investigations in this field did not give
any indication of such a serious obstacle in the general theory of Hilbert
C*-modules. L.~G.~Brown obtained his examples from the theory of different
kinds of multipliers of C*-algebras without identity. Furthermore, making
use of the results of the Ph.D.~thesis of Nien-Tsu Shen \cite{Shen} he proved
the following: For a Banach C*-module $\cal M$ over a C*-algebra $A$ carrying
two $A$-valued inner products $\langle .,. \rangle_1$, $\langle
.,. \rangle_2$ which induce equivalent to the given one norms on $\cal M$
the appropriate C*-algebras of ''compact'' bounded $A$-linear operators on
the Hilbert $A$-modules $\{ {\cal M}, \langle .,. \rangle_1 \}$ and
$\{ {\cal M}, \langle .,. \rangle_2 \}$ are $*$-isomorphic if
and only if there exists a bounded $A$-linear isomorphism $S$ of these two
Hilbert $A$-modules satisfying $\langle .,. \rangle_2 \equiv
\langle S(.),S(.) \rangle_1$, cf.~\cite[Thm.~4.2, Prop.~4.4]{Brown} together
with \cite[Prop.~2.3]{Frank:93}, (\cite{BMS}). By definition, the set of
''compact'' operators ${\rm K}_A({\cal M})$ on a Hilbert $A$-module $\{
{\cal M}, \langle .,. \rangle \}$ is defined as the norm-closure of the set
${\rm K}_A^0({\cal M})$ of all finite linear combinations of the operators
    \[
    \{ \theta_{x,y} \, : \, \theta_{x,y}(z) = \langle z,x
    \rangle y \; \, {\rm for} \: {\rm every} \; \, x,y,z \in \cal M \}.
    \]
It is a C*-subalgebra and a two-sided ideal of ${\rm End}_A^*(\cal M)$, the set
of all adjointable bounded $A$-linear operators on $\{ {\cal M}, \langle .,.
\rangle \}$, what is the multiplier C*-algebra of ${\rm K}_A({\cal M})$ by
\cite[Thm.~1]{Kas}. Note, that in difference to the well-known situation for
Hilbert spaces, the properties of an operator to be ''compact'' or to possess
an adjoint depend heavily on the choice of the $A$-valued inner product on
$\cal M$. These properties are not invariant even up to isomorphic Hilbert
structures on $\cal M$, in general, cf.~\cite{Frank:93}. We make the
convention that operators $T$ which are ''compact''/adjointable with respect
to some $A$-valued inner product $\langle .,. \rangle_i$ will be marked
$T^{(i)}$ to note where this property arises from. The same will be done for
sets of such operators.

\noindent
In 1994 E.~C.~Lance showed that two Hilbert C*-modules are isomorphic as
Hilbert C*-modules if and only if they are isometrically isomorphic as
Banach C*-modules (\cite{Lance:94}) opening the geometrical background
of this functional-analytical problem and extending a central result for
C*-algebras: C*-algebras are isometrically multiplicatively isomorphic if
and only if they are $*$-isomorphic, \cite[Thm.~7, Lemma 8]{Kad}.

\smallskip \noindent
At the contrary,
non-isomorphic Hilbert structures on a given Hilbert $A$-module $\cal M$
over a C*-algebra $A$ can not appear at all if $\cal M$ is {\it self-dual},
i.~e.~every bounded module map $r: {\cal M} \rightarrow A$ is of the form
$\langle . , a_r \rangle$ for some element $a_r \in \cal M$
(cf.~\cite[Prop.~2.2,Cor.~2.3]{Fr1}), or if $A$ is unital and $\cal M$ is
{\it countably generated}, i.~e.~there exists a countably set of generators
inside $\cal M$ such that the set of all finite $A$-linear combinations of
generators is norm-dense in $\cal M$  (cf.~\cite[Cor.~4.8, Thm.~4.9]{Brown}
together with \cite[Cor.~1.1.25]{JT} and \cite[Prop.~2.3]{Frank:93}).

\medskip \noindent
Now, we come to the goal of the present paper: Whether for a Banach C*-module
$\cal M$ over a C*-algebra $A$ carrying two $A$-valued inner products
$\langle .,. \rangle_1$, $\langle .,. \rangle_2$ which induce equivalent to
the given one norms on $\cal M$ the appropriate C*-algebras
${\rm End}_A^{(1,*)}({\cal M})$ and ${\rm End}_A^{(2,*)}({\cal M})$ of all
adjointable bounded $A$-linear operators on $\cal M$ are $*$-isomorphic, or
not? This question is non-trivial since even non-$*$-isomorphic non-unital
C*-algebras can possess a common multiplier C*-algebra:
For example, on the closed interval $[0,2] \subset {\bf R}$ consider
the C*-algebra of all continuous functions vanishing at zero together
with the C*-algebra of all continuous function vanishing at one.
They are non-$*$-isomorphic, but the multiplier C*-algebra C([0,2]) of them
consisting of all continuous functions on [0,2] is the same in both cases.
That is, additional arguments are needed to describe the relation between
the multiplier C*-algebras of non-$*$-isomorphic C*-algebras of ''compact''
operators on some Banach C*-modules carrying non-isomorphic C*-valued inner
products. One quickly realizes that the techniques of multiplier theory are
not suitable to shed some more light on this general situation. One has to
turn back to C*-theory and to the properties of $*$-isomorphisms, as
well as to the theory of Hilbert C*-modules.


\medskip \noindent
{\bf Theorem:} $\:$ {\it
Let $A$ be a C*-algebra and $\cal M$ be a Banach $A$-module carrying two
$A$-valued inner products $\langle .,. \rangle_1$, $\langle .,. \rangle_2$
which induce equivalent to the given one norms. Then the following conditions
are equivalent:
  \newcounter{cou001}
  \begin{list}{(\roman{cou001})}{\usecounter{cou001}}
  \item The Hilbert $A$-modules $\{ {\cal M}, \langle .,. \rangle_1 \}$
        and $\{ {\cal M}, \langle .,. \rangle_2 \}$ are isomorphic as Hilbert
        C*-modules.
  \item The Hilbert $A$-modules $\{ {\cal M}, \langle .,. \rangle_1 \}$
        and $\{ {\cal M}, \langle .,. \rangle_2 \}$ are isometrically
        isomorphic as Banach $A$-modules.
  \item The C*-algebras ${\rm K}_A^{(1)}({\cal M})$ and
        ${\rm K}_A^{(2)}({\cal M})$ of all ''compact'' bounded $A$-linear
        ope\-rators on both these Hilbert C*-modules, respectively, are
        $*$-isomorphic.
  \item The unital C*-algebras ${\rm End}_A^{(1,*)}({\cal M})$ and
        ${\rm End}_A^{(2,*)}({\cal M})$ of all adjointable bounded $A$-linear
        operators on both these Hilbert C*-modules, respectively, are
        $*$-isomorphic.
  \end{list}
Further equivalent conditions in terms of positive invertible quasi-multipliers
of ${\rm K}_A^{(1)}({\cal M})$ can be found in {\rm \cite{Frank:93}}.
}

\medskip \noindent
{\sc Proof.} The equivalence of (i) and (ii) was shown by E.~C.~Lance
\cite{Lance:94}, and the equivalence of (i) and (iii) turns out from
a result for C*-algebras of L.~G.~Brown \cite[Thm.~4.2, Prop.~4.4]{Brown}
in combination with \cite[Prop.~2.3]{Frank:93}. Referring to G.~G.~Kasparov
\cite[Thm.~1]{Kas} the implication (iii)$\to$(iv) yields naturally.

\noindent
Now, suppose the unital C*-algebras ${\rm End}_A^{(1,*)}({\cal M})$
and ${\rm End}_A^{(2,*)}({\cal M})$ are $*$-isomorphic. Denote this
$*$-isomorphism by $\omega$. One quickly checks that the formula
\[
x \in {\cal M} \rightarrow \langle x,x \rangle_{1,Op.} = \theta^{(1)}_{x,x}
\in {\rm K}^{(1)}_A({\cal M})
\]
defines a ${\rm K}^{(1)}_A(\cal M)$-valued inner product on the Hilbert
$A$-module $\cal M$ regarding it as a right ${\rm K}^{(1)}_A(\cal M)$-module.
Moreover, the set $\{ K(x) : x \in {\cal M}, K \in {\rm K}^{(1)}_A(\cal M)
\}$ is norm-dense inside $\cal M$ since the limit equality
\[
x=\|.\|_{\cal M}-\lim_{n \to \infty}
(\theta_{x,x}^{(1)}(\theta_{x,x}^{(1)}+n^{-1})^{-1})(x)
\]
holds for every $x \in \cal M$.

\noindent
As a first step we consider the intersection of the two C*-subalgebras and
two-sided ideals $\omega({\rm K}^{(1)}_A(\cal M))$ and ${\rm K}^{(2)}_A
(\cal M)$ inside the unital C*-algebra ${\rm End}_A^{(2,*)}({\cal M})$.
The intersection of them is a C*-subalgebra and two-sided ideal of
${\rm End}_A^{(2,*)}({\cal M})$ again. It contains the operators
\[
\theta^{(2)}_{x,y} \cdot \omega(\theta^{(1)}_{z,t})  =
                 \theta^{(2)}_{\omega(\theta^{(1)}_{z,t})^*(x),y}
                  =   \theta^{(2)}_{\omega(\theta^{(1)}_{t,z})(x),y}
\]
for every $x,y,z,t \in \cal M$. Since the set of all finite linear
combinations of special opera\-tors $\{ \theta^{(1)}_{z,t} \, : \, z,t \in \cal
M \}$ is norm-dense inside ${\rm K}^{(1)}_A(\cal M)$ by definition the
intersection of $\omega({\rm K}^{(1)}_A(\cal M))$ and ${\rm K}^{(2)}_A(\cal M)$
contains the set
\[
\{ \theta^{(2)}_{\omega(K^{(1)})(x),y} : K^{(1)} \in {\rm K}^{(1)}_A
({\cal M}) \: , \: x,y \in {\cal M} \}  \, .
\]
Because of the limit equality
\begin{eqnarray*}
x & = & \|.\|_{\cal M}-\lim_{n \to \infty}
\omega(\theta_{x,x}^{(1)}(\theta_{x,x}^{(1)}+n^{-1})^{-1})(x) \\
  & = & \|.\|_{\cal M}-\lim_{n \to \infty}
\omega(\theta_{x,x}^{(1)})\omega((\theta_{x,x}{(1)}+n^{-1})^{-1})(x)
\end{eqnarray*}
the set $\{ \omega(K^{(1)})(x) : K^{(1)} \in {\rm K}^{(1)}_A({\cal M}),
x \in {\cal M} \}$ is norm-dense inside $\cal M$.
Consequently, the intersection of $\omega({\rm K}^{(1)}_A(\cal M))$ and
${\rm K}^{(2)}_A(\cal M)$ inside the unital C*-algebra ${\rm End}_A^{(2,*)}
({\cal M})$ contains the set of ''compact'' operators $\{ \theta^{(2)}_{x,y}
\, : \, x,y \in \cal M \}$ generating one of the intersecting sets,
${\rm K}^{(2)}_A(\cal M)$, completely, and the inclusion relation
${\rm K}^{(2)}_A({\cal M}) \subseteq \omega({\rm K}^{(1)}_A({\cal M}))$
holds.

\noindent
Secondly, by the symmetry of the situation and of the arguments the inclusion
relation ${\rm K}^{(1)}_A({\cal M}) \subseteq \omega^{-1}({\rm K}^{(2)}_A
({\cal M}))$ holds, too, inside the unital C*-algebra ${\rm End}_A^{(1,*)}
({\cal M})$. Both inclusions together prove that $\omega$ realizes a
$*$-isomorphism of the C*-algebras ${\rm K}^{(1)}_A({\cal M})$ and
${\rm K}^{(2)}_A({\cal M})$  automatically, what implies (iii) and hence, (i).
$\: \bullet$

\medskip \noindent
Whether the $*$-isomorphism of the C*-algebras of ''compact'' bounded
$A$-linear operators of two different Hilbert $A$-modules $\cal M$ and
$\cal N$ over some C*-algebras $A$ implies their isomorphism as Hilbert
C*-modules, or not? The answer is negative, even in the quite well-behaved
cases. Counterexamples appear because of nontrivial $K_0$-groups of $A$,
for instance. Let $A$ be the hyperfinite type ${\rm II}_1$ W*-factor.
Set ${\cal M}=A$ and ${\cal N}=A^2$ with the usual $A$-valued inner products.
Both these Hilbert $A$-modules are self-dual and finitely generated.
Obviously, ${\rm K}_A({\cal M})$ and ${\rm K}_A({\cal N})$ are $*$-isomorphic
to $A$ as C*-algebras. Nevertheless, $\cal M$ and $\cal N$ are not isomorphic
as Banach $A$-modules because of the non-existence of non-unitary isometries
for the identity caused by the existence of a faithful trace functional on $A$.
The $K_0$-group of $A$ equals ${\bf R}$, i.~e., it is non-trivial, and
$A \cong A \otimes {\rm M}_2({\bf C})$.

\noindent
In general, one could search for some special unital C*-algebra $A$ with
non-trivial $K_0$-group, a natural number $n \geq 1$ and two projections $p,q
\in {\rm M}_n(A)$ such that for every $N \geq n$ the finitely generated
Hilbert $A$-modules $A^Np$ and $A^Nq$ are non-isomorphic
(i.~e., $[p] \not= [q] \in K_0(A)$), but the C*-algebras $p {\rm M}_n(A) p$
and $q {\rm M}_n(A) q$ are $*$-isomorphic.

\medskip \noindent
Closing, we pose the problem whether for a Banach C*-module $\cal M$ over a
C*-algebra $A$ carrying two $A$-valued inner products $\langle .,. \rangle_1$,
$\langle .,. \rangle_2$ which induce equivalent to the given one norms on
$\cal M$ the appropriate Banach algebras of all (not necessarily adjointable)
bounded $A$-linear operators on $\cal M$ are {\it isometrically}
multiplicatively isomorphic, or not, especially in the case of non-isomorphic
Hilbert structures. Those properties of all these kinds of operator algebras
which are preserved switching from one $A$-valued inner product on $\cal M$
to another have to be investigated in the future extending results for the
''compact'' case of \cite{BMS,Frank:93}.

\newpage

{\small


Fed.~Rep.~Germany

Universit\"at Leipzig

FB Mathematik/Informatik

Mathematisches Institut

Augustusplatz 10

D-04109 Leipzig

frank@mathematik.uni-leipzig.d400.de

}
\end{document}